\documentclass[12pt]{article}
\usepackage[a3paper]{geometry}
\usepackage[numbers]{natbib}
\usepackage{graphicx}
\usepackage{color, multirow,amsfonts,amsmath}
\usepackage{url}
\usepackage{rotating}
\begin{document}
\title{Bivariate compound Poisson risk processes with shocks\\({\it{Working version}})
\\ Pavlina K. Jordanova
\\{\small{\it Faculty of Mathematics and Informatics, Konstantin Preslavsky University of Shumen, \\115 "Universitetska" str., 9712 Shumen, Bulgaria.
\\ Corresponding author: pavlina\_kj@abv.bg}}
\\ Evelina Veleva
\\{\small{\it{Department of Applied mathematics and Statistics, "Angel Kanchev" University of Ruse, Bulgaria.}}}
\\ Kosto Mitov
\\{\small{\it{Department of Medical Physics, Biophysics, Pre-clinical and Clinical Sciences, Medical University - Pleven, Bulgaria.}}}}
\maketitle

\begin{abstract}
Contemporary insurance theory is concentrated on models with different types of polices and shock events may influence the payments on some of them.  Jordanova \cite{JordanovaMattex2018} considered a model where a shock event contributes to the total claim amount with one and the same value of the claim sizes to different types of polices. Jordanova and Veleva \cite{Jordanova2021} went a step closer to real life situations and allowed a shock event to cause different claim sizes to different types of polices. In that paper the counting process is assumed to be Multinomial. Here it is replaced with different independent homogeneous Poison processes. The bivariate claim counting process is expressed in two different ways. Its marginals and conditional distributions are totaly described. The mean square regression of these processes is computed. The Laplace-Stieltjes transforms and numerical characteristics of the total claim amount processes are obtained.  The risk reserve process and the probabilities of ruin in infinite time are discussed. The risk reserve just before the ruin and the deficit (or the severity) at ruin are thoroughly investigated in case when the initial capital is zero. Their means, probability mass functions and probability generating functions are obtained.

Although the model is constructed by a multivariate counting processes, along the paper it is shown that the total claim amount process is stochastically equivalent to a univariate compound Poisson process. These allows us to reduce the considered risk model to a Cramer-Lundberg risk model, to use the corresponding results and to make the conclusions for the new model. Analogous results can be obtained for more types of polices and more types of shock events.

The results are applied in case when the claim sizes are exponentially distributed.

Stochastically equivalent models could be analogously constructed in queuing theory.
\end{abstract}

\section{DESCRIPTION OF THE MODEL AND ITS RELATIONS WITH PREVIOUSLY INVESTIGATED MODELS}

The origin of modeling of random processes with common shocks is usually related with the names of Marshal and Olkin \cite{Marshal&Olkin1967,Marshal&Olkin1988} or Mardia \cite{Mardia1970}. They construct different bivariate models with dependent coordinates, where the dependence is caused by a joint summand, or a joint term in the maxima or the minima. The results could be easily transferred to random processes with independent increments. Cossette and Marceau \cite{CossetteMarceau2000} apply such results in risk theory. Jordanova \cite{JordanovaMattex2018} uses merging of compound Poisson processes and investigates a particular case of the model, defined in this work. The difference is that, in that work, in a fixed point of time, if there is a shock event it affects different types of contracts via one and the same claim size. In the insurance model investigated in Jordanova and Veleva \cite{Jordanova2021} the form of the dependence is improved. The claim sizes caused by one and the same shock event on different types of polices can be different. However, in that paper, the counting processes of claim arrivals of different types of polices are described via a Multinomial processes.  Here this Multinomial process is replaced with different independent homogeneous Poisson processes (HPPs) and the form of the dependence is the same. More precisely, let $N_0$, $N_1$, $N_2$ be mutually independent HPPs with parameters correspondingly $\lambda_0 > 0$, $\lambda_1 > 0$ and $\lambda_2 > 0$. By assumption $N_0(0) = N_1(0) = N_2(0) = 0$.

The  counting processes $M_1$ and $M_2$ of both types of claims are in generally dependent HPP with common shocks, and for all $ t \geq 0$
\begin{equation}\label{M}
M_1(t) = N_1(t) +  N_0(t), \quad M_2(t) = N_2(t) +  N_0(t),
\end{equation}
Therefore, $M_1$ is a HPP with parameter $\lambda_1 + \lambda_0$ and $M_2$ is a HPP with parameter $\lambda_2 + \lambda_0$.

The total claim amount process $S = \{S(t):  t \geq 0\}$ is such that for all $t \geq 0$,
 \begin{eqnarray}
\label{S}   S(t) &=& S_1(t) + S_2(t),\\
\label{S1S2} S_1(t) &=& \sum_{i = 0}^{N_1(t)} Y_{1i} + \sum_{i = 0}^{N_0(t)} Y_{3i},\quad Y_{10} = 0,\,\, Y_{30} = 0\\
\nonumber S_2(t) &=& \sum_{i = 0}^{N_2(t)} Y_{2i} + \sum_{i = 0}^{N_0(t)} Y_{4i},\quad Y_{20} = 0,\,\, Y_{40} = 0.
 \end{eqnarray}

The sequences $\{Y_{ki}\}_{i = 1}^\infty$ and $\{N_r(s), s\geq 0\}$,  $k = 1, 2, 3, 4$, $r = 0, 1, 2$ are mutually independent.
The processes $S_1 = \{S_1(t), \,  t \geq 0\}$ and $S_2 = \{S_2(t), \, t \geq 0\}$ are in generally dependent.

The definitions of $S$, $S_1$ and $S_2$ mean that  we allow the claims to be of two mutually exclusive different types and shock events (like car-crashes) may cause possibly different payments in both polices. The  process  $N_0 = \{N_0(t), \, t \geq 0\}$ describes the counts of shock events.
In (\ref{S1S2}) for any fixed $k = 1, 2, 3, 4$ the random variables (r.vs.) $Y_{k1}, Y_{k2}, ...$ are independent identically distributed (i.i.d.)  with cumulative distribution function (c.d.f.) $\mathbb{F}_k$. The sequences $\{Y_{1i}\}_{i=1}^\infty$, $\{Y_{2i}\}_{i=1}^\infty$ and $\{(Y_{3i},Y_{4i})\}_{i=1}^\infty$ are mutually independent. Without lost of generality the almost sure strictly positive r.v. $Y_{ki}$ describes the $i$-th claim amount, ($i \in \mathbb{N}$) of a certain type $k$, $k \in \{1, 2, 3, 4\}$. If the distribution of the claim sizes has mass at zero, the parameters of the corresponding counting HPPs, participating in $S_1$ and $S_2$, defined in (\ref{S1S2}), can be changed, and the model can be reduced to the considered here one.

A natural example of an insurance interpretation of this model is as follows. For $k = 1$, $Y_{1i}$ is the part of the claim size which is due to the customer for $i$-th standard insurance policy of type $1$ exclusively the common shocks. Analogously for $k = 2$. The r.v. $Y_{3i} + Y_{4i}$ models the part of the claim size which have to be paid to the customer due to some common shock. Jordanova and Veleva \cite{Jordanova2021} give and example when the random elements with $k = 1$ describe health insurances, $k = 2$ is for car insurance. Then, $\sum_{i = 0}^{N_1(t)} Y_{1i}$ is the total claim amount up to time $t$, which is due to customers for having some healthy problems, not related with car-crashes. The r.v. $\sum_{i = 0}^{N_0(t)} Y_{3i}$ models the total claim amount up to time $t$, which is due to healthy problems, caused by car-crashes. The r.v. $\sum_{i = 0}^{N_2(t)} Y_{2i}$ is for the total claim amount up to time $t$, which is due to some car problems which are not related with car-crashes, for example  car-thefts. And finally, the r.v.
$\sum_{i = 1}^{N_0(t)} Y_{4i}$ is the total claim amount up to time $t$, which is due to some car problems which are related with car-crashes.
Therefore, the random process $S_1$ is the total claim amount processes for the insurance events caused by "healthy problems", $S_2$ is the amount paid to customers for events caused by "problems with cars".

As usually the risk reserve process $R_u = \{R_u(t): t \geq 0\}$ is defined via the equality
\begin{equation}\label{R}
R_u(t) = u + ct - S(t), \quad t \geq 0,
\end{equation}
where $u \geq 0$ is the initial capital and $c > 0$ is the premium income rate per unit time.

The time of ruin when the initial capital is $u$ will be denoted by
\begin{eqnarray}
\label{timeofruin}
\tau(u)&:=& \inf\{t > 0: R_u(t) < 0\},\quad  \inf\, \emptyset = \infty,
\end{eqnarray}

The probability for ruin in "infinite horizon"  and initial capital $u$ is
\begin{equation}\label{probforruin}
\psi(u) := \mathbb{P}(\tau(u) < \infty).
\end{equation}
and the corresponding survival probability is $\delta(u) := 1 - \psi(u)$.

All random elements discussed here are assumed to be measurable with respect to one the same probability space  $(\Omega, \mathcal{F}, \mathbb{P})$ with natural filtration born by the considered processes. The letter $\mathbb{G}$ is a notation for probability generating function (p.g.f.). Analogously, the letter $\mathbb{L}$ means Laplace-Stieltjes transform (LST).

\bigskip
\section{THE BIVARIATE COUNTING PROCESS}

Poisson and compound Poisson processes are very well investigated in the scientific literature. See for example Kingman \cite{Kingman1993}, or Sundt and Vernic \cite{SundtVernic}. In multivariate case we have two types of compound Poisson processes. The I type means "with equal number of summands". The type "II" is its opposite. Let us now derive the compound Poisson  presentation of the first type of the bivariate counting process $(M_1, M_2) := \{(M_1(t), M_2(t)): t \geq 0\}$, defined in (\ref{M}). In order to formulate our results we will define a random vector $(I_1, I_2, I_0)$ with probability mass function (p.m.f.)
\begin{equation}\label{pmf1}
\mathbb{P}(I_1 = 1, I_2 = 0, I_0 = 0) = \frac{\lambda_1}{\lambda}, \,\,  \mathbb{P}(I_1 = 0, I_2 = 1, I_0 = 0) = \frac{\lambda_2}{\lambda}, \,\,\mathbb{P}(I_1 = 0, I_2 = 0, I_0 = 1) = \frac{\lambda_0}{\lambda},
\end{equation}
and zero otherwise. Here $\lambda_0 > 0$, $\lambda_1 > 0$, $\lambda_2 > 0$ and  $\lambda: = \lambda_0 + \lambda_1 + \lambda_2$. By assumption the random vectors $(I_{11}, I_{21}, I_{01}), (I_{12}, I_{22}, I_{02}),  ...$ are i.i.d. with p.m.f. (\ref{pmf1}).

Denote by, $N = \{N(t), \, t \geq 0\}$ a HPP with $N(0) = 0$, and parameter $\lambda$, and by
\begin{equation}
A_1(t) := \sum_{i = 0}^{N(t)} (I_{1i} + I_{0i}),\quad A_2(t) := \sum_{i = 0}^{N(t)} (I_{2i} + I_{0i}), \quad t \geq 0. \label{A}
\end{equation}
Here $I_{10} = I_{20} = I_{00} = 0$ and  $N$ and $(I_{11}, I_{21}, I_{01}), (I_{12}, I_{22}, I_{02}),  ...$ are independent.  In our example the random process $A_1$ counts the number of the insurance events caused by healthy problems. The random process $A_2$ counts the number of the insurance events caused by problems with cars.

In the next theorem we show that the bivariate counting process
$(M_1, M_2) := \{(M_1(t), M_2(t)): t \geq 0\}$ is a particular case of a multivariate homogeneous in time compound Poisson process with equal number of summands, i.e. of type I, in the sense of Sundt and Vernic \cite{SundtVernic,Vernic2018}.

\medskip

{\bf Theorem 1.} The bivariate counting processes $(M_1, M_2)$ and $(A_1, A_2)$ coincide in the sense of their finite dimensional distributions.

{\bf Proof:} The processes $(M_1, M_2)$ and $(A_1, A_2)$ have homogeneous and independent additive increments and they start from the coordinate beginning. Therefore, in order to prove their stochastic equivalence it is enough to prove equality of their univariate time intersections. Due to the uniqueness of the correspondence between the probability laws and their p.g.fs. it is enough to derive equality between the p.g.fs.

The definition of p.g.fs., (\ref{M}), the multiplicative property of p.g.f. and  the well-known formula for the p.g.f. of the Poisson distribution entail that for all $t \geq 0$,
\begin{eqnarray*}
&\mathbb{E}&\left[z_1^{M_1(t)}z_2^{M_2(t)}\right] = \mathbb{E}\left[z_1^{N_1(t) +  N_0(t)}z_2^{N_2(t) +  N_0(t)}\right]  = \mathbb{E}\left[z_1^{N_1(t)}z_2^{N_2(t)} (z_1z_2)^{N_0(t)} \right]   =  \mathbb{E}\left[z_1^{N_1(t)} \right]\mathbb{E}\left[z_2^{N_2(t)} \right] \mathbb{E}\left[(z_1z_2)^{N_0(t)} \right]\\
&=& \exp\{-\lambda_1t(1-z_1)\}\exp\{-\lambda_2t(1-z_2)\}\exp\{-\lambda_0t(1-z_1z_2)\} = \exp(-\lambda t+\lambda_1tz_1+\lambda_2tz_2+\lambda_0tz_1z_2).
\end{eqnarray*}

By the definition (\ref{A}), the formula for double expectations and the definition of Poisson distribution
\begin{eqnarray*}
\mathbb{E}\left[z_1^{A_1(t)}z_2^{A_2(t)}\right] &=&  \mathbb{E}\left[z_1^{\sum_{i = 0}^{N(t)} (I_{1i} + I_{0i})} z_2^{\sum_{i = 0}^{N(t)} (I_{2i} + I_{0i})} \right] = \sum_{n=0}^\infty \mathbb{E}\left[z_1^{\sum_{i = 0}^{n} (I_{1i} + I_{0i})} z_2^{\sum_{i = 0}^{n} (I_{2i} + I_{0i})} \right]\mathbb{P}(N(t) = n)\\
& =& \sum_{n=0}^\infty \mathbb{E}\left[z_1^{\sum_{i = 0}^{n} (I_{1i} + I_{0i})} z_2^{\sum_{i = 0}^{n} (I_{2i} + I_{0i})} \right]\frac{(\lambda t)^n}{n!}e^{-\lambda t}.
\end{eqnarray*}

By the multiplicative property of p.g.fs. and the fact that the random vectors $(I_{11}, I_{21}, I_{01}), (I_{12}, I_{22}, I_{02}),  ...$ are i.i.d. we have
\begin{equation}
\mathbb{E}\left[z_1^{A_1(t)}z_2^{A_2(t)}\right] = \sum_{n=0}^\infty \left[\mathbb{E}\left(z_1^{I_{11} + I_{01}} z_2^{I_{21} + I_{01}} \right)\right]^n \frac{(\lambda t)^n}{n!}e^{-\lambda t}.
\end{equation}

Now, by using (\ref{pmf1}) we obtain the p.g.f. of $\mathbb{E}\left[z_1^{I_{11} + I_{01}} z_2^{I_{21} + I_{01}} \right]$.
Therefore, the fact that $N$ is a HPP with parameter $\lambda$ and the Maclaurin series leads us to
\begin{eqnarray*}
\mathbb{E}\left[z_1^{A_1(t)}z_2^{A_2(t)}\right] &=& \sum_{n=0}^\infty \left[ \frac{\lambda_1}{\lambda}z_1 + \frac{\lambda_2}{\lambda}z_2 + \frac{\lambda_0}{\lambda} z_1z_2\right]^n \frac{(\lambda t)^n}{n!}e^{-\lambda t} = \exp(-\lambda t+\lambda_1tz_1+\lambda_2tz_2+\lambda_0tz_1z_2).
\end{eqnarray*}
which is exactly $\mathbb{E}\left[z_1^{M_1(t)}z_2^{M_2(t)}\right]$. The uniqueness of the correspondence between p.g.fs. and the probability distribution completes the proof. \hfill $\Box$

{\it Note 1.} For all $t \geq 0$ the distribution of the time intersections $(M_1(t), M_2(t))$ of the bivariate counting process $(M_1, M_2)$, defined in (\ref{M}) are determined via the following p.g.fs.,
$$\mathbb{G}_{M_1(t),M_2(t)}(z_1, z_2) = \mathbb{E}(z_1^{M_1(t)}z_2^{M_2(t)}) = \exp(-\lambda t+\lambda_1tz_1+\lambda_2tz_2+\lambda_0tz_1z_2).$$

They are called bivariate Poisson distributions, and for any fixed $t \geq 0$, they are very well investigated in the scientific literature. See for example Kocherlakota and Kocherlakota \cite{Kocherlakota1992} who show that for $m_1 , m_2 \in \{0, 1, ...\}$,
$$\mathbb{P}(M_1(t) = m_1, M_2(t) = m_2) = e^{-\lambda t}\sum_{i = 0}^{min(m_1, m_2)} \frac{(\lambda_0t)^{i}(\lambda_1t)^{m_1-i}(\lambda_2t)^{m_2-i}}{i!(m_1-i)!(m_2-i)!},$$
and this probability mass function (p.m.f.) is equal to $0$ otherwise.

Thus, for $(r,s) = (1,2)$ or $(r,s) = (2,1)$, $m_s \in \mathbb{N}$ and $m_r \in \mathbb{N}$
$$\mathbb{P}(M_r(t) = m_r|M_s(t) = m_s)= e^{-(\lambda_0+\lambda_r)t}m_s!\sum_{i = 0}^{min(m_1, m_2)} \frac{(\lambda_0t)^{i}(\lambda_rt)^{m_r-i}(\lambda_s t)^{-i}}{i!(m_1-i)!(m_2-i)!},$$
and this p.m.f. is equal to $0$ otherwise.
$$\mathbb{E}(M_i(t)) = \mathbb{D}(M_i(t)) = t(\lambda_i + \lambda_0), \,i = 1, 2, \quad \mathbb{E}(M_1(t)M_2(t)) = \lambda_0 t + (\lambda \lambda_0 +  \lambda_1 \lambda_2) t^2, \quad cov(M_1(t),M_2(t)) = t\lambda_0.$$

 Their correlation does not depend on $t$ and is
 $$cor(M_1(t),M_2(t)) = \frac{\lambda_0}{\sqrt{(\lambda_1 + \lambda_0)(\lambda_2 + \lambda_0)}}.$$

{\bf Theorem 2.} For $(r,s) = (1,2)$ or $(r,s) = (2,1)$, $m_s \in \mathbb{N}$ and $m_r \in \mathbb{N}$ the mean square regression of the bivariate couunting processes $(M_1, M_2)$ is
$$\mathbb{E}(M_s(t)|M_r(t) = m) = \frac{m\lambda_0}{\lambda_r+\lambda_0} +\lambda_st, \quad  m = 0, 1, ...$$
{\bf Proof:} Without lost of generality let us assume that  $(r,s) = (1,2)$. Then,
$$\frac{\partial^m\mathbb{G}_{M_1(t),M_2(t)}(z_1, z_2)}{\partial z_1^m}  = \exp(-\lambda t+\lambda_1tz_1+\lambda_2tz_2+\lambda_0tz_1z_2)(\lambda_1t+\lambda_0tz_2)^m,$$
$$\frac{\partial^{m+1}\mathbb{G}_{M_1(t),M_2(t)}(z_1, z_2)}{\partial z_1^{m} \partial z_2}  = \exp(-\lambda t+\lambda_1tz_1+\lambda_2tz_2+\lambda_0tz_1z_2)\left[(\lambda_1t+\lambda_0tz_2)^m(\lambda_2t+\lambda_0tz_1)+ m(\lambda_1t+\lambda_0tz_2)^{m-1}\lambda_0t\right],$$
$$\mathbb{E}(M_2(t)|M_1(t) = m) = \frac{\left.\frac{\partial^{m+1}\mathbb{G}_{M_1(t),M_2(t)}(z_1, z_2)}{\partial z_1^{m} \partial z_2} \right|_{(z_1, z_2) = (0, 1)}}{\left.\frac{\partial^m\mathbb{G}_{M_1(t),M_2(t)}(z_1, z_2)}{\partial z_1^m}\right|_{(z_1, z_2) = (0, 1)}} = \frac{\exp(-\lambda t+\lambda_2t)\left[(\lambda_1t+\lambda_0t)^m\lambda_2t - m(\lambda_1t+\lambda_0t)^{m-1}\lambda_0t\right]}{\exp(\lambda t+\lambda_2t)(\lambda_1t+\lambda_0t)^m}$$
$$\mathbb{E}(M_2(t)|M_1(t) = m) = \frac{m\lambda_0}{\lambda_1+\lambda_0} + \lambda_2t.$$
\hfill $\Box$

{\it Note  2.} The process $A = M_1 + M_2 = N_1 + N_2 + 2N_0$ is a Markov chain with $Q$-matrix
$$Q = \left(\begin{array}{ccccccc}
                                                                                              -\lambda & \lambda_1 + \lambda_2 & \lambda_0 & 0 & 0 & 0 &... \\
                                                                                              0 & -\lambda & \lambda_1 + \lambda_2 & \lambda_0 & 0 & 0 &... \\
                                                                                              0 & 0 & -\lambda & \lambda_1 + \lambda_2 & \lambda_0 & 0 & ... \\
                                                                                              ... & ... & ... &... & ... & ... & ...\\
                                                                                            \end{array}
                                                                                          \right).$$
It is a univariate compound Poisson process. For all $t \geq 0$, its time intersections satisfy the equalities
$$A(t) := M_1(t) + M_2(t) = \sum_{i = 0}^{N(t)} (I_{1i} + I_{2i} + 2I_{0i}).$$
It is determined via the distribution of its time intersections which have the following p.g.fs.,
$$\mathbb{E}(z^{A(t)}) = \exp[-t(\lambda-(\lambda_1+\lambda_2)z-\lambda_0z^2)], \quad z \geq 0.$$
 It is easy to see  that as far as $\mathbb{E}(I_{1i} + I_{2i} + 2I_{0i}) = 1+\frac{\lambda_0}{\lambda}$ and $\mathbb{D}(I_{1i} + I_{2i} + 2I_{0i}) = \frac{\lambda_0}{\lambda}\left(1 - \frac{\lambda_0}{\lambda}\right)$, thus
$$\mathbb{E}(A(t)) = t(\lambda_1 + \lambda_2 + 2\lambda_0), \quad \mathbb{D}(A(t)) = (\lambda_1 + \lambda_2 + 4\lambda_0)t.$$

\bigskip
\section{STOCHASTICALLY EQUIVALENT PRESENTATIONS OF THE TOTAL CLAIM AMOUNT PROCESS}

The process $(S_1, S_2)$ has a compound Poisson of type I presentation.

\medskip

{\bf Theorem 3.} The bivariate total claim amount process $(S_1, S_2)$, described in (\ref{S1S2}) is a bivariate compound Poisson process of type I. It is stochastically equivalent (in the sense of finite dimensional distributions) to the process $(S_3, S_4)$, where
\begin{equation}\label{S4}
S_3(t) = \sum_{n = 0}^{N(t)}(I_{1n}Y_{1n} + I_{0n}Y_{3n}), \quad S_4(t) = \sum_{n = 0}^{N(t)}(I_{2n}Y_{2n} + I_{0n}Y_{4n}), \quad t \geq 0.
\end{equation}
Here $Y_{k0} = 0$, $k = 1, 2, 3, 4$, $(I_{11}, I_{21}, I_{01})$ is a vector with distribution (\ref{pmf1}), and $N$ is a HPP with parameter $\lambda = \lambda_1 + \lambda_2 + \lambda_0$. The random process $N$, and the sequences $\{Y_{1i}\}_{i=1}^\infty$, $\{Y_{2i}\}_{i=1}^\infty$, $\{(Y_{3i}, Y_{4i})\}_{i=1}^\infty$, and  $\{(I_{1i}, I_{2i}, I_{0i})\}_{i=1}^\infty$ are mutually independent.

{\bf Proof:}  Consider $t \geq 0$. The definition of Laplace-Stieltjes transforms (LSTs), (\ref{S1S2}),  the multiplicative property of LSTs,  and the formula for the LSTs of compound Poisson distribution entail
\begin{eqnarray*}
 && \mathbb{E}e^{-z_1S_1(t) - z_2S_2(t)} = \mathbb{E}e^{-z_1\sum_{i = 0}^{N_1(t)} Y_{1i} - z_1 \sum_{i = 0}^{N_0(t)} Y_{3i} - z_2 \sum_{i = 0}^{N_0(t)} Y_{4i}  - z_2 \sum_{i = 0}^{N_2(t)} Y_{2i}}\\
&=& \mathbb{E}e^{-z_1 \sum_{i = 0}^{N_1(t)} Y_{1i}}\mathbb{E}e^{-\sum_{i = 0}^{N_0(t)}(z_1Y_{31}+ z_2Y_{41})}\mathbb{E}e^{- z_2 \sum_{i = 0}^{N_2(t)} Y_{2i}}\\
&=& \exp\{-\lambda_1 t(1-\mathbb{E}e^{-z_1Y_{11}})\}\exp\{-\lambda_0 t(1-\mathbb{E}e^{-(z_1Y_{31}+ z_2Y_{41})})\exp\{-\lambda_2 t(1-\mathbb{E}e^{-z_2Y_{21}})\}\\
&=&  \exp\{-\lambda_1 t(1-\mathbb{E}e^{-z_1Y_{11}})-\lambda_0 t(1-\mathbb{E}e^{-(z_1Y_{31}+ z_2Y_{41})})-\lambda_2 t(1-\mathbb{E}e^{-z_2Y_{21}})\}\\
&=&  \exp\left\{-\lambda t\left(1-\frac{\lambda_1}{\lambda}\mathbb{E}e^{-z_1Y_{11}}-\frac{\lambda_0}{\lambda}\mathbb{E}e^{-(z_1Y_{31}+ z_2Y_{41})}-\frac{\lambda_2}{\lambda}\mathbb{E}e^{-z_2Y_{21}}\right)\right\}\\
&=&  \exp\left\{-\lambda t\left(1-\mathbb{E}e^{-(z_1I_{11}Y_{11} + z_1I_{01}Y_{31} + z_2I_{21}Y_{21} + z_2I_{01}Y_{41})}\right)\right\}\\
&=& \mathbb{E} \exp\left\{-\sum_{i = 0}^{N(t)}[z_1(I_{1i}Y_{1i} + I_{0i}Y_{3i})+z_2(I_{2i}Y_{2i} + I_{0i}Y_{4i})]\right\}\\
&=&\mathbb{E} \exp\left\{-z_1\sum_{i = 0}^{N(t)}(I_{1i}Y_{1i} + I_{0i}Y_{3i}) - z_2\sum_{i = 0}^{N(t)}(I_{2i}Y_{2i} + I_{0i}Y_{4i})\right\} = \mathbb{E}e^{-z_1S_3(t) - z_2S_4(t)}.
\end{eqnarray*}
In the last equality we have used the definitions (\ref{S4}).  The uniqueness of the correspondence between the probability laws and their LSTs entails the equality in distribution of the univariate time intersections of $(S_1, S_2)$, described in (\ref{S1S2}), and those of the process $(S_3, S_4)$. These processes have homogeneous and independent additive increments and they start from the coordinate beginning. Therefore, the last entails their stochastic equivalence in the sense of finite dimensional distributions.   \hfill $\Box$

{\bf Corollary 1.} For all $t \geq 0$ the distribution of the time intersections $(S_1(t), S_2(t))$, $t = 1, 2, ...$ of the bivariate process $(S_1, S_2)$ defined in (\ref{S1S2}) are determined via the following LSTs,
$$\mathbb{E}e^{-z_1S_1(t) - z_2S_2(t)} = \exp\left\{-\lambda\left(1-\frac{\lambda_1}{\lambda}\mathbb{E}e^{-z_1Y_{11}}-\frac{\lambda_0}{\lambda}\mathbb{E}e^{-(z_1Y_{31}+ z_2Y_{41})}-\frac{\lambda_2}{\lambda}\mathbb{E}e^{-z_2Y_{21}}\right)\right\}.$$
\begin{itemize}
\item[i)] If  $\mathbb{E}Y_{i1} < \infty$, $i = 1, 2, 3, 4$, then, $\mathbb{E}(S_1(t)) = t(\lambda_1\mathbb{E}Y_{11} + \lambda_0\mathbb{E}Y_{31})$, $\mathbb{E}(S_2(t)) = t(\lambda_2\mathbb{E}Y_{21} + \lambda_0\mathbb{E}Y_{41});$
\item[ii)] If $\mathbb{E}(Y_{i1}^2) < \infty$, $i = 1, 2, 3, 4$, then,
$$\mathbb{D}(S_1(t)) = t\left[\lambda_1\mathbb{E}(Y_{11}^2) + \lambda_0\mathbb{E}(Y_{31}^2)\right], \quad \mathbb{D}(S_2(t)) =  t\left[\lambda_2\mathbb{E}(Y_{21}^2) + \lambda_0\mathbb{E}(Y_{41}^2)\right];$$
\item[iii)]
$\mathbb{E}(S_1(t)S_2(t)) = t\lambda_0\mathbb{E}(Y_{31}Y_{41})+ t^2(\lambda_1\mathbb{E}Y_{11} + \lambda_0\mathbb{E}Y_{31})(\lambda_2\mathbb{E}Y_{21} + \lambda_0\mathbb{E}Y_{41})$;
\item[iv)]
$cov(S_1(t), S_2(t)) =  t\lambda_0\mathbb{E}(Y_{31}Y_{41});$
\item[v)] The correlation $cor(S_1(t),S_2(t))$  does not depend on $t$. More precisely
$$cor(S_1(t),S_2(t)) = \frac{\lambda_0\mathbb{E}(Y_{31}Y_{41})}{\sqrt{\left[\lambda_1\mathbb{E}(Y_{11}^2) + \lambda_0\mathbb{E}(Y_{31}^2)\right]\left[\lambda_2\mathbb{E}(Y_{21}^2) + \lambda_0\mathbb{E}(Y_{41}^2)\right]}}.$$
\end{itemize}

{\bf Corollary 2.} The time intersections
$$S(t) := \sum_{i = 0}^{N(t)} \left[I_{1i}Y_{1i} + I_{0i}(Y_{3i}+Y_{4i}) + I_{2i}Y_{2i}\right], \quad t \geq 0,$$
of the compound Poisson process $S = \{S(t), \,t\geq 0\}$, possess the following properties:
\begin{itemize}
\item[i)] their distributions are determined via the following Laplace-Stieltjes transforms
$$\mathbb{E}(e^{-zS(t)}) = \exp\left\{-\lambda t\left[1 - \frac{\lambda_1}{\lambda}\mathbb{E}e^{-zY_{11}} - \frac{\lambda_0}{\lambda}\mathbb{E}e^{-z(Y_{31}+Y_{41})} - \frac{\lambda_2}{\lambda} \mathbb{E}e^{-zY_{21}}\right]\right\}, \quad t \geq 0.$$
\item[ii)] $\mathbb{E}(S(t)) = t[\lambda_1\mathbb{E}Y_{11} + \lambda_0(\mathbb{E}Y_{31} + \mathbb{E}Y_{41})+ \lambda_2\mathbb{E}Y_{21}]$;
\item[iii)] $\mathbb{D}(S(t)) = t\left\{\lambda_1\mathbb{E}(Y_{11}^2) + \lambda_0\mathbb{E}[(Y_{31}+ Y_{41})^2] + \lambda_2\mathbb{E}(Y_{21}^2) \right\}.$
\end{itemize}

\bigskip
\section{STOCHASTICALLY EQUIVALENT RISK MODELS}
In this section, first it is shown that the risk reserve process $R_u$, defined in (\ref{R}), is a classical compound Poisson risk process. Then, by using the results about the Cramer-Lundberg risk model, the corresponding characteristics of $R_u$  and probabilities for ruin in infinite time are obtained. More precisely,
Corollary 2 entails that the risk process, described in (\ref{R}) is stochastically equivalent to the process
\begin{equation}\label{Rt}
R_u(t) = u + ct - \sum_{i = 0}^{N(t)} Y_i, \quad t \geq 0,
\end{equation}
where $Y_0 = 0$, and $Y_i: = I_{1i}Y_{1i} + I_{0i}(Y_{3i} + Y_{4i}) + I_{2i}Y_{2i}$. Let us denote by  $\mu : = \mathbb{E} Y_1 = \frac{\lambda_1}{\lambda}\mathbb{E}Y_{11} + \frac{\lambda_2}{\lambda}\mathbb{E}Y_{21} + \frac{\lambda_0}{\lambda}(\mathbb{E}Y_{31} + \mathbb{E}Y_{41})$. By assumption $\mu  < \infty$. Then,
$$\mathbb{E}[R_u(t)] = u + ct - \lambda t \mu = u + ct - t(\lambda_1\mathbb{E}Y_{11} + \lambda_2\mathbb{E}Y_{21} + \lambda_0(\mathbb{E}Y_{31} + \mathbb{E}Y_{41})).$$

By using these notations we obtain the well-known formula for the safety loading in the Cramer-Lundberg risk model,
\begin{equation}\label{safetyloading}
\rho:= \lim_{t \to \infty} \frac{\mathbb{E}[R_u(t)]}{\mathbb{E}[S(t)]} = \frac{c}{\lambda\mu}-1 = \frac{c}{\lambda_1\mathbb{E}Y_{11} + \lambda_2\mathbb{E}Y_{21} + \lambda_0(\mathbb{E}Y_{31} + \mathbb{E}Y_{41})}-1.
\end{equation}
The corresponding net profit condition is
\begin{equation}\label{npc}
\rho > 0 \iff c > \lambda\mu \iff c > \lambda_1\mathbb{E}Y_{11} + \lambda_2\mathbb{E}Y_{21} + \lambda_0(\mathbb{E}Y_{31} + \mathbb{E}Y_{41}).
\end{equation}
For the general results about the Cramer-Lundberg risk model see for example Rolski et al. \cite{Rolski1998},  Asmussen and Albrecher \cite{AsmussenAlbrecher2010}, Grandell \cite{Grandell1991} or  Gerber \cite{Gerber1979} among others.

If the net profit condition (\ref{npc}) is satisfied, then the Cramer-Lundberg presentation (\ref{Rt}) of the risk process $R_u$ allows us to conclude that it possesses the following properties.
$$\psi(0) = \frac{1}{1 + \rho} = \frac{\lambda \mu}{c} = \frac{\lambda_1}{c}\mathbb{E}Y_{11} + \frac{\lambda_2}{c}\mathbb{E}Y_{21} + \frac{\lambda_0}{c}(\mathbb{E}Y_{31} + \mathbb{E}Y_{41}),$$
$$\delta(0) = \frac{\rho}{1 + \rho} = \frac{c - \lambda_1\mathbb{E}Y_{11} - \lambda_2\mathbb{E}Y_{21} - \lambda_0(\mathbb{E}Y_{31} + \mathbb{E}Y_{41})}{c}$$

If we denote by $Z_0 = 0$, by $Z_k = \sum_{i=1}^k (Y_i - cX_i)$, $k = 1, 2, ...$, and by $M = \sup_{i \in \mathbb{N}}Z_i$, then
$\delta(u) = \mathbb{P}(M \leq u)$, and
 \begin{eqnarray}
 \label{delta} \delta(u) &=& \delta(0) + \frac{\lambda}{c} \int_0^u \delta(u-y)(1-F_{Y_1}(y))dy\\
\nonumber &=& \delta(0) + \frac{\lambda_1}{c}\int_0^u \delta(u-y)(1-F_{Y_{11}}(y))dy + \frac{\lambda_0}{c}\int_0^u \delta(u-y)(1-F_{Y_{31} + Y_{41}}(y))dy\\
\nonumber & +& \frac{\lambda_2}{c}\int_0^u \delta(u-y)(1-F_{Y_{21}}(y))dy.\\
\label{deltaprim} \delta'(u) &=& \frac{\lambda}{c}\left(\delta(u)-\int_0^u \delta(u-y)dF_{Y_1}(y)\right)\\
\nonumber &=& \frac{\lambda}{c}\delta(u)- \frac{\lambda_1}{c}\int_0^u \delta(u-y)dF_{Y_{11}}(y)- \frac{\lambda_0}{c}\int_0^u \delta(u-y)dF_{Y_{31} + Y_{41}}(y) - \frac{\lambda_2}{c}\int_0^u \delta(u-y)dF_{Y_{21}}(y).
\end{eqnarray}

Now, we have the compound geometric form, for the Laplace-Stieltes transform of $M$. For $s > 0$,
\begin{eqnarray}
\nonumber   \mathbb{E}(e^{-sM}) &=& \int_0^\infty e^{-su}\delta'(u)du  + \delta(0) =  \frac{\delta(0)}{1-(1-\delta(0))\left(\frac{1 -  \mathbb{E}(e^{-sY_{1}})}{\mu s}\right)}\\
\label{LST} &=& \frac{(c-\lambda_1\mathbb{E}Y_{11} - \lambda_2\mathbb{E}Y_{21} - \lambda_0(\mathbb{E}Y_{31} + \mathbb{E}Y_{41}))s}{cs - \lambda  + \lambda_1\mathbb{E}(e^{-sY_{11}}) + \lambda_0 \mathbb{E}(e^{-s(Y_{31} + Y_{41})}) + \lambda_2 \mathbb{E}(e^{-sY_{21}})}
 \end{eqnarray}

Given that $\tau(0)$ is finite, let us denote by $\eta$ the absolute value of the deficit at ruin. This r.v. has an integrated tail distribution
 $$\mathbb{P}(\eta < x) = \mathbb{P}(-R_0(\tau(0)^+) \leq x|\tau(0)< \infty) = \frac{\int_0^x(\lambda_1 \mathbb{P}(Y_{11} \geq y) + \lambda_2 \mathbb{P}(Y_{21} \geq y)  + \lambda_0  \mathbb{P}(Y_{31} + Y_{41} \geq y) )dy}{\lambda_1\mathbb{E}Y_{11} + \lambda_2\mathbb{E}Y_{21} + \lambda_0(\mathbb{E}Y_{31} + \mathbb{E}Y_{41})}=: F_I(x)$$
 with mean value
$$\mathbb{E}(-R_0(\tau(0)^+)|\tau(0) < \infty) = \frac{\lambda_1\mathbb{E}(Y_{11}^2) + \lambda_2\mathbb{E}(Y_{21}^2) + \lambda_0 \mathbb{E}(Y_{31} + Y_{41})^2}{2\left(\lambda_1\mathbb{E}Y_{11} + \lambda_2\mathbb{E}Y_{21} + \lambda_0(\mathbb{E}Y_{31} + \mathbb{E}Y_{41})\right)}$$ and the Laplace-Stieltes transform
 $$ \mathbb{E}(e^{-s\eta}) = \frac{1 -  \mathbb{E}(e^{-sY_{1}})}{\mu s} = \frac{\lambda - \lambda_1\mathbb{E}(e^{-sY_{11}}) - \lambda_0 \mathbb{E}(e^{-s(Y_{31} + Y_{41})}) - \lambda_2 \mathbb{E}(e^{-sY_{21}})}{(\lambda_1\mathbb{E}Y_{11} + \lambda_2\mathbb{E}Y_{21} + \lambda_0(\mathbb{E}Y_{31} + \mathbb{E}Y_{41}))s}, \quad s > 0.$$

 As far as the distribution of $\eta$ coincides with the integrated tail distribution of the mixture $Y_1$ of the distributions of $Y_{11}$, $Y_{21}$ and $Y_{31} + Y_{41}$, it is a mixture of the corresponding integrated tail distributions of $Y_{11}$, $Y_{21}$ and $Y_{31} + Y_{41}$. The corresponding weights are
 $$p_i := \frac{\lambda_i\mathbb{E}Y_{i1}}{\lambda_1\mathbb{E}Y_{11} + \lambda_2\mathbb{E}Y_{21} + \lambda_0(\mathbb{E}Y_{31} + \mathbb{E}Y_{41})}, \quad i = 1, 2, \quad  p_0 := \frac{\lambda_0(\mathbb{E}Y_{31}+\mathbb{E}Y_{41})}{\lambda_1\mathbb{E}Y_{11} + \lambda_2\mathbb{E}Y_{21} + \lambda_0(\mathbb{E}Y_{31} + \mathbb{E}Y_{41})}.$$

 {\bf Theorem 4.} If $\mathbb{E}Y_{i1} < \infty$, $i = 1, 2, 3, 4$, and the random vector $(J_{11}, J_{21}, J_{01})$ has the distribution
  \begin{equation}\label{pi}
  \mathbb{P}(J_{11} = 1, J_{21} = 0, J_{01} = 0) = p_1, \quad \mathbb{P}(J_{11} = 0, J_{21} = 1, J_{01} = 0) = p_2, \quad \mathbb{P}(J_{11} = 0, J_{21} = 0, J_{01} = 1) = p_0,
  \end{equation}
    then the integrated tail distribution of  $Y_1: = I_{11}Y_{11} + I_{01}(Y_{31} + Y_{41}) + I_{21}Y_{21}$ is the mixture with c.d.f.
    $$\frac{1}{\mathbb{E}Y_1}\int_0^y \mathbb{P}(Y_1 \geq x)dx = \frac{p_1}{\mathbb{E}Y_{11}}\int_0^\infty \mathbb{P}(Y_{11} \geq x)dx + \frac{p_0}{\mathbb{E}Y_{31}+\mathbb{E}Y_{41}}\int_0^y \mathbb{P}(Y_{31} + Y_{41} \geq x)dx + \frac{p_2}{\mathbb{E}Y_{21}}\int_0^\infty \mathbb{P}(Y_{21} \geq x)dx.$$

     The last means that if the r.vs. $\theta_1, \theta_2, \theta_0$ have the integrated tail distribution correspondingly of $Y_{11}$, $Y_{21}$, and $Y_{31} + Y_{41}$, and if the r.vs. $\theta_1, \theta_2, \theta_0$ are independent of $(J_{11}, J_{21}, J_{01})$, then
     $$\eta \stackrel{d}{=} J_{11}\theta_1 + J_{21}\theta_2 + J_{01}\theta_0.$$

    The proof follows by the total probability formula.

 \medskip

 In this way   $$M \stackrel{d}{=} \sum_{i = 0}^{\xi} \eta_i,$$
where the r.vs. $\eta, \eta_1, \eta_2, ...$ are independent and identically distributed, and they are independent on $\xi$. The r.v. $\xi$ is Geometrically distributed over the natural numbers ($\mathbb{P}(\xi = k) = q^{k-1}p$, $k \in \mathbb{N}$), with parameter
 $$p =  \frac{c - \lambda_1\mathbb{E}Y_{11} - \lambda_2\mathbb{E}Y_{21} - \lambda_0(\mathbb{E}Y_{31} + \mathbb{E}Y_{41})}{c}.$$
Now, we have obtained the Beekmans convolution series \cite{Beekman1968} or Pollaczek-Khinchin formula \cite{Khinchine1932}
 $$\delta(u) = \delta(0) \sum_{i=0}^\infty (1-\delta(0))^i F_\eta^{i*}(u), \quad u\geq 0.$$

The expected time of ruin, given that the ruin occur and the initial capital is $u = 0$ is
$$\mathbb{E}(\tau(0)|\tau(0) < \infty) = \frac{\lambda_1\mathbb{E}(Y_{11}^2) + \lambda_2\mathbb{E}(Y_{21}^2) + \lambda_0 \mathbb{E}(Y_{31} + \mathbb{E}Y_{41})^2}{2\left(\lambda_1\mathbb{E}Y_{11} + \lambda_2\mathbb{E}Y_{21} + \lambda_0(\mathbb{E}Y_{31} + \mathbb{E}Y_{41})\right)\left(c - \lambda_1\mathbb{E}Y_{11} + \lambda_2\mathbb{E}Y_{21} + \lambda_0(\mathbb{E}Y_{31} + \mathbb{E}Y_{41})\right)}.$$

The joint distribution of the deficit at ruin and $R(\tau(0)^-)$ - the risk surplus just before the ruin, when $u = 0$ is
\begin{eqnarray*}
\mathbb{P}(-R_0(\tau(0)^+) > x, R(\tau(0)^-) > y|\tau(0)< \infty) &=& \frac{\int_{x+y}^\infty(\lambda_1 \mathbb{P}(Y_{11} \geq y) + \lambda_2 \mathbb{P}(Y_{21} \geq y)  + \lambda_0  \mathbb{P}(Y_{31} + Y_{41} \geq y) )dy}{\lambda_1\mathbb{E}Y_{11} + \lambda_2\mathbb{E}Y_{21} + \lambda_0(\mathbb{E}Y_{31} + \mathbb{E}Y_{41})}.
\end{eqnarray*}

 In case when the initial capital is positive, when we replace c.d.f. of $Y_1$ in the well known formulae for $G(u,y):=\mathbb{P}(-R(\tau(u)+)\leq y, \tau(u)< \infty)$, see for example Gerber and Shiu (1997) \cite{GerberandShiu1997} or Klugman, Panjer and Willmot (2012) \cite{KPW2012},
$$\frac{\partial}{\partial y}G(u,y) = \frac{\lambda}{c}G(u,y) - \frac{\lambda_1}{c}\left[\int_0^u G(u-x,y)dF_{Y_{11}}(y) + F_{Y_{11}}(u+y) - F_{Y_{11}}(u)\right]$$
$$ - \frac{\lambda_2}{c}\left[\int_0^u G(u-x,y)dF_{Y_{21}}(y) + F_{Y_{21}}(u+y) - F_{Y_{21}}(u)\right] - \frac{\lambda_0}{c}\left[\int_0^u G(u-x,y)dF_{Y_{31} + Y_{41}}(y) + F_{Y_{31} + Y_{41}}(u+y) - F_{Y_{31} + Y_{41}}(u)\right].$$

$$G(u,y) = \frac{\lambda_1}{c}\left[\int_0^u G(u-x,y)\bar{F}_{Y_{11}}(x)dx +\int_u^{u+y}\bar{F}_{Y_{11}}(x)dx\right]+\frac{\lambda_2}{c}\left[\int_0^u G(u-x,y)\bar{F}_{Y_{21}}(x)dx +\int_u^{u+y}\bar{F}_{Y_{21}}(x)dx\right]$$
$$+\frac{\lambda_0}{c}\left[\int_0^u G(u-x,y)\bar{F}_{Y_{31} + Y_{41}}(x)dx +\int_u^{u+y}\bar{F}_{Y_{31} + Y_{41}}(x)dx\right].$$

By the main properties of the Cramer-Lundberg risk model, which could be found in any textbook on Risk theory, we have that if the Cramer-Lundberg exponent $\varepsilon > 0$ exists, which in this case is the same, as if a positive solution of the equation
$$\lambda_1\mathbb{E}e^{\varepsilon Y_{11}}+\lambda_2\mathbb{E}e^{\varepsilon Y_{21}} + \lambda_0\mathbb{E}e^{\varepsilon(Y_{31}+Y_{41})} - \lambda = \frac{c \lambda \varepsilon}{\lambda_1\mathbb{E}Y_{11} + \lambda_2\mathbb{E}Y_{21} + \lambda_0(\mathbb{E}Y_{31} + \mathbb{E}Y_{41})}$$
 exists, then
$\psi(u) \leq e^{-\varepsilon u}$.

If additionally
$$\int_0^\infty ye^{\varepsilon y}(\lambda_1\mathbb{P}(Y_{11} \geq y) + \lambda_2 \mathbb{P}(Y_{21} \geq y)  + \lambda_0  \mathbb{P}(Y_{31} + Y_{41} \geq y) )dy < \infty,$$
then
$$\lim_{u \to \infty} e^{\varepsilon u}\psi(u) = \frac{c - \lambda_1\mathbb{E}Y_{11} - \lambda_2\mathbb{E}Y_{21} - \lambda_0(\mathbb{E}Y_{31} + \mathbb{E}Y_{41})}{\varepsilon(\int_0^\infty ye^{\varepsilon y}(\lambda_1\mathbb{P}(Y_{11} \geq y) + \lambda_2 \mathbb{P}(Y_{21} \geq y)  + \lambda_0  \mathbb{P}(Y_{31} + Y_{41} \geq y) )dy}.$$

If the distributions of $Y_{11}$, $Y_{21}$, $Y_{31}$, and $Y_{41}$ are all subexponential and $\mathbb{E}Y_{11}$, $\mathbb{E}Y_{21}$, $\mathbb{E}Y_{31}$, and $\mathbb{E}Y_{41}$ are all finite, then
$$\lim_{u \to \infty} \frac{\psi(u)}{\int_0^u(\lambda_1 \mathbb{P}(Y_{11} \geq y) + \lambda_2 \mathbb{P}(Y_{21} \geq y)  + \lambda_0  \mathbb{P}(Y_{31} + Y_{41} \geq y) )dy} = \frac{1}{c - \lambda_1\mathbb{E}Y_{11} - \lambda_2\mathbb{E}Y_{21} - \lambda_0(\mathbb{E}Y_{31} + \mathbb{E}Y_{41})}$$

\section{Exponential claim sizes}

Let us now suppose that $Y_{k1}$, $k = 1, 2, 3, 4$ are independent and exponentially distributed with means correspondingly $\nu_k$, $k = 1, 2, 3, 4$. Then,
\begin{equation}\label{safetyloadingExp}
\rho = \frac{c}{\lambda_1\nu_1 + \lambda_2\nu_2 + \lambda_0(\nu_3 + \nu_4)}-1,
\end{equation}
$$\mu : = \frac{\lambda_1\nu_1}{\lambda} + \frac{\lambda_2\nu_2}{\lambda} + \frac{\lambda_0}{\lambda}(\nu_3 + \nu_4), $$
\begin{equation}\label{piExp}
p_i := \frac{\lambda_i\nu_i}{\lambda_1\nu_1 + \lambda_2\nu_2 + \lambda_0(\nu_3 + \nu_4)}, \quad i = 1, 2, \quad  p_0 := \frac{\lambda_0(\nu_3+\nu_4)}{\lambda_1\nu_1 + \lambda_2\nu_2 + \lambda_0(\nu_3 + \nu_4)}.
\end{equation}
The corresponding net profit condition is
\begin{equation}\label{npcExp}
c > \lambda_1\nu_1 + \lambda_2\nu_2 + \lambda_0(\nu_3 + \nu_4).
\end{equation}

If the net profit condition (\ref{npc}) is satisfied,
$$\psi(0) = \frac{\lambda_1}{c}\nu_1 + \frac{\lambda_2}{c}\nu_2 + \frac{\lambda_0}{c}(\nu_3 + \nu_4),$$
$$\delta(0) = 1 - \frac{\lambda_1}{c}\nu_1 + \frac{\lambda_2}{c}\nu_2 + \frac{\lambda_0}{c}(\nu_3 + \nu_4)$$
 \begin{eqnarray}
\nonumber  \delta(u) &=& \delta(0) + \frac{\lambda_1}{c}\int_0^u \delta(u-y)e^{-\frac{y}{\nu_1}}dy + \frac{\lambda_0}{c}\frac{\nu_3}{\nu_3-\nu_4}\int_0^u \delta(u-y)e^{-\frac{y}{\nu_3}}dy -\frac{\lambda_0}{c}\frac{\nu_4}{\nu_3-\nu_4}\int_0^u \delta(u-y)e^{-\frac{y}{\nu_4}}dy\\
\label{deltaExp} & +& \frac{\lambda_2}{c}\int_0^u \delta(u-y)e^{-\frac{y}{\nu_2}}dy.\\
\nonumber \delta'(u) &=& \frac{\lambda}{c}\delta(u)+ \frac{\lambda_1}{c}\int_0^u \delta(u-y)de^{-\frac{y}{\nu_1}}+\frac{\nu_3}{\nu_3-\nu_4} \frac{\lambda_0}{c}\int_0^u \delta(u-y)de^{-\frac{y}{\nu_3}}- \frac{\lambda_0}{c}\frac{\nu_4}{\nu_3-\nu_4}\int_0^u \delta(u-y)de^{-\frac{y}{\nu_4}}\\
\label{deltaprimExp} &+& \frac{\lambda_2}{c}\int_0^u \delta(u-y)de^{-\frac{y}{\nu_2}}
\end{eqnarray}

Now, we have the compound geometric form, for the Laplace-Stieltes transform (LST) of $M$. For $s > 0$,
\begin{eqnarray}
 \label{LSTExp}  \mathbb{E}(e^{-sM}) &=& \int_0^\infty e^{-su}\delta'(u)du  + \delta(0) = \frac{(c-\lambda_1\nu_1 - \lambda_2\nu_2 - \lambda_0(\nu_3 + \nu_4))s}{cs - \lambda  + \frac{\lambda_1}{1 + s\nu_1} + \frac{\lambda_0}{(1+s\nu_3)(1+s\nu_4)} + \frac{\lambda_2}{1+\nu_2}}.
 \end{eqnarray}

As far as the integrated tail distribution of the Exponential distribution is again Exponential distribution with the same parameter, the integrated tail distribution of $Y_1$ has a LST
\begin{eqnarray}\label{etaLSTExp}
\mathbb{E}e^{-s\eta} &=& p_1 \frac{1-\mathbb{E}e^{-sY_{11}}}{s\mathbb{E}Y_{11}} + p_2 \frac{1-\mathbb{E}e^{-sY_{21}}}{s\mathbb{E}Y_{21}} + p_0 \frac{1-\mathbb{E}e^{-s(Y_{31}+Y_{41})}}{s\mathbb{E}(Y_{31}+Y_{41})}\\
&=& p_1 \frac{1}{1+\nu_1s} + p_2 \frac{1}{1+\nu_1s} + p_0 \frac{1-\frac{1}{1+\nu_3s}\frac{1}{1+\nu_4s}}{s(\nu_3+\nu_4)}\\
&=& p_1 \frac{1}{1+\nu_1s} + p_2 \frac{1}{1+\nu_1s} + p_0 \frac{\nu_3s + \nu_4s + \nu_2\nu_4s^2}{s(\nu_3+\nu_4)(1-\nu_3s)(1-\nu_4s)}.
\end{eqnarray}
Now we need to observe only that
$$\frac{\nu_3s + \nu_4s + \nu_2\nu_4s^2}{s(\nu_3+\nu_4)(1-\nu_3s)(1-\nu_4s)}$$ is the LST of the integrated tail distribution of $Y_{31} + Y_{41}$. Let us suppose that a r.v. $\theta_0$ has such probability law, and $\theta_0$ is independent of $Y_{11}$ and $Y_{21}$. Then
     $$\eta \stackrel{d}{=} J_{11}Y_{11} + J_{21}Y_{21} + J_{01}\theta_0.$$

      The mean value of $\eta$ is,
$$\mathbb{E}\eta = \mathbb{E}(-R_0(\tau(0)^+)|\tau(0) < \infty) = \frac{\lambda_1\nu_1^2 + \lambda_2\nu_2^2 + \lambda_0(\nu_3^2 + \nu_3\nu_4 + \nu_4^2)}{\lambda_1\nu_1 + \lambda_2\nu_2 + \lambda_0(\nu_3 + \nu_4))}.$$

In this case $\eta$ it is a mixture of two Exponential and the integrated tail distribution of a Hypoexponential distribution. Moreover,
\begin{eqnarray}\label{theta0}
\mathbb{P}(\theta_0 \leq x) &=& \frac{\int_0^x \mathbb{P}(Y_{31}+Y_{41} \geq z)dz}{\mathbb{E}Y_{31}+\mathbb{E}Y_{41} } = \frac{x - \int_0^x \int_0^z \mathbb{P}(Y_{31} < z-y)d\mathbb{P}(Y_{41} < y)dz}{\mathbb{E}Y_{31}+\mathbb{E}Y_{41} }\\
\nonumber  &=& \frac{x -  z\int_0^z \mathbb{P}(Y_{31} < z-y)d\mathbb{P}(Y_{41} < y)|_{z =0}^x + \int_0^x z \int_0^z P_{Y_{31}}(z-y)d\mathbb{P}(Y_{41} < y)dz}{\mathbb{E}Y_{31}+\mathbb{E}Y_{41} }\\
\nonumber  &=& \frac{x -  x\int_0^x \mathbb{P}(Y_{31} < x-y)d\mathbb{P}(Y_{41} < y) + \int_0^x z \int_0^z P_{Y_{31}}(z-y)d\mathbb{P}(Y_{41} < y)dz}{\mathbb{E}Y_{31}+\mathbb{E}Y_{41} }.
\end{eqnarray}
If $\nu_3 \not= \nu_4$, then for $x > 0$,
\begin{eqnarray}\label{theta0Exp}
\mathbb{P}(\theta_0 \leq x) &=& \frac{x -  \frac{x}{\nu_4}\int_0^x (1-e^{\frac{-(x-y)}{\nu_3}})e^{-\frac{y}{\nu_4}}dy + \frac{1}{\nu_3\nu_4}\int_0^x z \int_0^z e^{-\frac{z-y}{\nu_3}}e^{-\frac{y}{\nu_4}}dydz}{\nu_3+\nu_4}\\
\nonumber &=& \frac{xe^{-\frac{x}{\nu_4}} +  \frac{xe^{-\frac{x}{\nu_3}}}{\nu_4}\int_0^x e^{y\left(\frac{1}{\nu_3}-\frac{1}{\nu_4}\right)}dy + \frac{1}{\nu_3\nu_4}\int_0^x ze^{-\frac{z}{\nu_3}} \int_0^z e^{y\left(\frac{1}{\nu_3}-\frac{1}{\nu_4}\right)}dydz}{\nu_3+\nu_4}\\
\nonumber &=& \frac{xe^{-\frac{x}{\nu_4}} +  \frac{xe^{-\frac{x}{\nu_3}}}{\nu_4\left(\frac{1}{\nu_3}-\frac{1}{\nu_4}\right)}\left(e^{x\left(\frac{1}{\nu_3}-\frac{1}{\nu_4}\right)}-1\right) + \frac{1}{\nu_3\nu_4\left(\frac{1}{\nu_3}-\frac{1}{\nu_4}\right)}\int_0^x ze^{-\frac{z}{\nu_3}} \left(e^{z\left(\frac{1}{\nu_3}-\frac{1}{\nu_4}\right)}-1\right)dz}{\nu_3+\nu_4}\\
\nonumber &=& \frac{xe^{-\frac{x}{\nu_4}} +  \frac{x}{\nu_4\left(\frac{1}{\nu_3}-\frac{1}{\nu_4}\right)}\left(e^{-\frac{x}{\nu_4}}-e^{-\frac{x}{\nu_3}}\right) + \frac{1}{\nu_4 -\nu_3}\int_0^x z \left(e^{-\frac{z}{\nu_4}}-e^{-\frac{z}{\nu_3}} \right)dz}{\nu_3+\nu_4}\\
\nonumber &=& \frac{\frac{x\nu_4e^{-\frac{x}{\nu_4}}}{\nu_4-\nu_3} - \frac{x\nu_3}{\nu_4-\nu_3}e^{-\frac{x}{\nu_3}} + \frac{1}{\nu_4 -\nu_3}\int_0^x z \left(e^{-\frac{z}{\nu_4}}-e^{-\frac{z}{\nu_3}} \right)dz}{\nu_3+\nu_4}\\
\nonumber &=& \frac{x\nu_4e^{-\frac{x}{\nu_4}} - x\nu_3e^{-\frac{x}{\nu_3}} + \int_0^x z e^{-\frac{z}{\nu_4}}dz-\int_0^x ze^{-\frac{z}{\nu_3}} dz}{\nu_4^2-\nu_3^2}\\
\nonumber &=& \frac{\nu_4^2(1- e^{-\frac{z}{\nu_4}}) - \nu_3^2(1- e^{-\frac{z}{\nu_3}})}{\nu_4^2-\nu_3^2} = 1 -  \frac{\nu_4^2e^{-\frac{z}{\nu_4}}}{\nu_4^2-\nu_3^2}  + \frac{\nu_3^2e^{-\frac{z}{\nu_3}}}{\nu_4^2-\nu_3^2},
\end{eqnarray}
and this c.d.f. is equal to $0$, otherwise. The corresponding probability density function is
$$P_{\theta_0}(x) = \left\{\begin{array}{ccc}
                             0 & , & x < 0 \\
                             \frac{\nu_4e^{-\frac{z}{\nu_4}}}{\nu_4^2-\nu_3^2} - \frac{\nu_3e^{-\frac{z}{\nu_3}}}{\nu_4^2-\nu_3^2} & , & x \geq 0
                           \end{array}
\right.$$

If $\nu_3 = \nu_4 = \nu_0$, then $\theta_0$ has the integrated tail distribution of $Gamma(2, \frac{1}{\nu_0})$ r.v. More precisely for $x > 0$,
\begin{eqnarray}\label{theta0Exp}
\mathbb{P}(\theta_0 \leq x) &=& \frac{1}{2\nu_0}\int_0^x\int_y^\infty \frac{1}{\nu_0^2}ze^{-\frac{z}{\nu_0}}dzdy \\
\nonumber &=& 1 - \frac{xe^{-\frac{x}{\nu_0}}}{2\nu_0} - e^{\frac{-x}{\nu_0}} =  \frac{1}{2}\left(1 - \frac{xe^{-\frac{x}{\nu_0}}}{\nu_0} - e^{\frac{-x}{\nu_0}}\right) + \frac{1}{2}\left(1 -  e^{\frac{-x}{\nu_0}}\right),
\end{eqnarray}
and this c.d.f. is equal to $0$, otherwise.
The last means that when $\nu_3 = \nu_4 = \nu_0$, $$\theta_0 \stackrel{d}{=} I_A \theta_3 + I_{\bar{A}}\theta_4,$$ where
\begin{equation}\label{A}
\mathbb{P}(I_A = 1, I_{\bar{A}} = 0) = \mathbb{P}(I_A = 0, I_{\bar{A}} = 1) = \frac{1}{2},
\end{equation}
 $\theta_3 \in Gamma(2, \frac{1}{\nu_0})$ with p.d.f.
$P_{\theta_3}(x)=  \frac{x}{\nu_0^2}e^{-\frac{x}{\nu_0}}$, when $x \geq 0$ and $P_{\theta_4}(x) = 0$, otherwise,
and $\theta_4 \in Exp\left(\frac{1}{\nu_0}\right)$ with p.d.f. $P_{\theta_4}(x) = \frac{1}{\nu_0}e^{-\frac{x}{\nu_0}}$, when $x \geq 0$ and $P_{\theta_4}(x) = 0$, otherwise.

Thus, when $\nu_3 \not= \nu_4$,
\begin{eqnarray}
 \label{etaExp}\mathbb{P}(\eta \leq x) &=& \mathbb{P}(-R_0(\tau(0)^+) \leq x|\tau(0)< \infty) = p_1 (1-e^{-\frac{x}{\nu_1}}) + p_2 (1-e^{-\frac{x}{\nu_2}}) + p_0 \mathbb{P}(\theta_0 \leq x)\\
  \nonumber &=& 1-\left[\frac{\lambda_1\nu_1}{\lambda_1\nu_1 + \lambda_2\nu_2+ \lambda_0(\nu_3 + \nu_4)}e^{-\frac{x}{\nu_1}} + \frac{\lambda_2\nu_2}{\lambda_1\nu_1 + \lambda_2\nu_2+ \lambda_0(\nu_3 + \nu_4)}e^{-\frac{x}{\nu_2}}\right. \\
& +&\left. \frac{\lambda_0(\nu_3 + \nu_4)}{\lambda_1\nu_1 + \lambda_2\nu_2+ \lambda_0(\nu_3 + \nu_4)}\left(\frac{\nu_3^2}{\nu_3^2-\nu_4^2}e^{-\frac{x}{\nu_3}}-\frac{ \nu_4^2}{\nu_3^2-\nu_4^2}e^{-\frac{x}{\nu_4}}\right)\right].
 \end{eqnarray}
This is equivalent to
\begin{equation}\label{MinExpCase}
M \stackrel{d}{=} \sum_{i = 0}^{\xi} (J_{1i}Y_{1i} + J_{2i}Y_{2i} + J_{0i}\theta_{0i})
\end{equation}
where the r.v. $\xi$ is Geometrically distributed over the natural numbers ($\mathbb{P}(\xi = k) = q^{k-1}p$, $k \in \mathbb{N}$), with parameter
 $$p =  1 - \frac{1}{c}(\lambda_1\nu_1 + \lambda_2\nu_2 + \lambda_0(\nu_3 + \nu_4)).$$
 $(J_{1i}, J_{2i}, J_{0i})$, $i = 1, 2, ...$ are i.i.d. random vectors with distribution defined in (\ref{pi}) with corresponding $p_1, p_2$, and $p_0$, defined in (\ref{piExp}). All these random elements are assumed to be independent.

If $\nu_3 = \nu_4 = \nu_0$, then, for $x > 0$,
\begin{eqnarray}
 \label{etaExpequal}\mathbb{P}(\eta \leq x) &=& \mathbb{P}(-R_0(\tau(0)^+) \leq x|\tau(0)< \infty) = p_1 (1-e^{-\frac{x}{\nu_1}}) + p_2 (1-e^{-\frac{x}{\nu_2}}) + p_0 \mathbb{P}(\theta_0 \leq x)\\
  \nonumber &=& 1-\left[\frac{\lambda_1\nu_1}{\lambda_1\nu_1 + \lambda_2\nu_2+ 2\lambda_0\nu_0}e^{-\frac{x}{\nu_1}} + \frac{\lambda_2\nu_2}{\lambda_1\nu_1 + \lambda_2\nu_2+ 2\lambda_0\nu_0}e^{-\frac{x}{\nu_2}} + \frac{2\lambda_0\nu_0}{\lambda_1\nu_1 + \lambda_2\nu_2+ 2\lambda_0\nu_0}\left(\frac{x}{2\nu_0} + 1\right)e^{\frac{-x}{\nu_0}}\right],
 \end{eqnarray}
 and $\mathbb{P}(\eta \leq x) =0$, otherwise.

  This is equivalent to
\begin{equation}\label{MinExpCase}
M \stackrel{d}{=} \sum_{i = 0}^{\xi} (J_{1i}Y_{1i} + J_{2i}Y_{2i} + J_{0i}(I_{Ai}\theta_{3i} + I_{\bar{A}i}\theta_{4i})
\end{equation}
where the r.v. $\xi$ is Geometrically distributed over the natural numbers, with parameter
 $$p =  1 - \frac{1}{c}(\lambda_1\nu_1 + \lambda_2\nu_2 + 2\lambda_0\nu_0).$$
 $(J_{1i}, J_{2i}, J_{0i})$, $i = 1, 2, ...$ are i.i.d. random vectors with distribution defined in (\ref{pi}) with corresponding $p_1, p_2$, and $p_0$, defined in (\ref{piExp}). For $i = 1, 2, ...$,  $(I_{Ai}, I_{\bar{A}i}) \stackrel{d}{=} (I_{A}, I_{\bar{A}})$. The last vector is defined in (\ref{A}). For $i = 1, 2, ...$,  $\theta_{3i} \in Gamma(2, \frac{1}{\nu_0})$,  and $\theta_{4i} \in Exp\left(\frac{1}{\nu_0}\right)$. All these random elements are assumed to be independent.

Formula (\ref{MinExpCase}) allows us to simulate independent observations on $M$ without generating all sample paths of the risk process and to plot the probability to survive, which is actually estimated via the empirical c.d.f. of $M$. Moreover, we can estimate also its derivative, which coincides with the probability density function of $M$.

\begin{figure}
\begin{minipage}[t]{0.5\linewidth}
    \includegraphics[scale=.7]{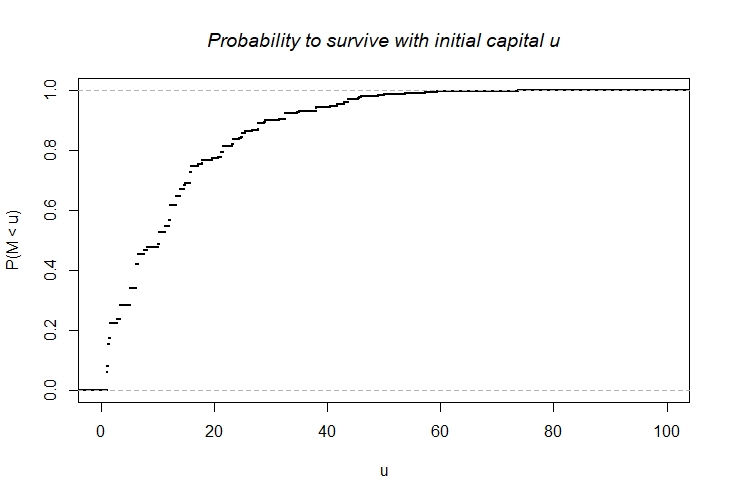}\vspace{-0.3cm}
\end{minipage}
\begin{minipage}[t]{0.49\linewidth}
    \includegraphics[scale=.7]{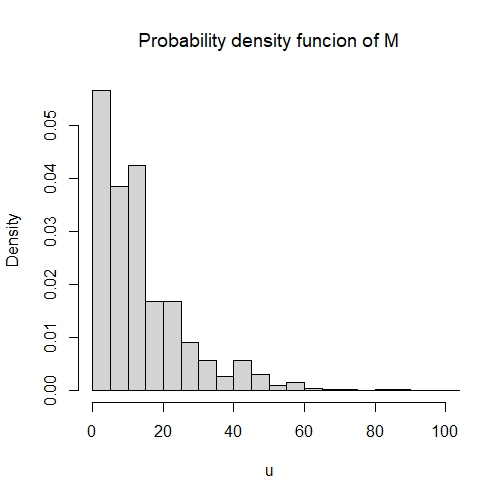}\vspace{-0.3cm}
\end{minipage}
\caption{\small  Estimators of the probability to survive (left) and its derivative (right) \label{fig:1}}
\end{figure}

Such simulation study was performed with sample size $10^7$ and parameters $\lambda_0 = 10$, $\lambda_1=11$, $\lambda_2 = 12$, $\nu_1 = 1$, $\nu_2 = 2$, $\nu_3 = \nu_4 = 3$, and $c = 97$. Then, $\rho \approx 0.021$, and $\mu = 95$. Both plots are given on Figure \ref{fig:1}.

\bigskip
\section{CONCLUSIONS}

Analogous conclusions could be done for mathematical models of insurance businesses with more types of claims, more types of polices, and more types of shock event. The idea is to reduce them to the classical compound Poisson risk model.

\bigskip

Although the model is constructed by a multivariate counting processes, along the paper it is shown that the total claim amount process is stochastically equivalent to a univariate compound Poisson process. These allows us to reduce the considered risk model to a Cramer-Lundberg risk model, to use the corresponding results and to make the conclusions for the new model.

\bigskip

These results give us a different approach for estimation of the probability of ultimate ruin and its derivative.
\bigskip

Analogous results can be applied in queuing theory.

\bigskip
\section{ACKNOWLEDGMENTS} The first author is partially supported by the project RD-08-144/01.03.2022 from the Scientific Research Fund in Konstantin Preslavsky University of Shumen, Bulgaria. Evelina Veleva thanks to the project No 2022-FNSE-04, financed by "Scientific research" Fund of Ruse University.

\bibliographystyle{aipnum-cp}%

\end{document}